\magnification=\magstep1
\hsize=16.5 true cm 
\vsize=23.6 true cm
\font\bff=cmbx10 scaled \magstep1

\font\bffg=cmbx10 scaled \magstep3

\parindent0cm
\overfullrule=0cm
%\nopagenumbers
%%%%%%%%%%%%%%%%%%%%%%%%%%%%%%%%
%%%%ABKšRZUNGEN im Textsatz %%%%
\def\cl{\centerline}           %
\def\rl{\rightline}            %
\def\bp{\bigskip}              %
\def\mp{\medskip}              %
\def\sp{\smallskip}            %
%%%%%%%%%%%%%%%%%%%%%%%%%%%%%%%%
%%%ZAHLENMENGEN%%%%%%%%%%%%%%%%%
\font\boldmas=msbm10           %
\def\Bbb#1{\hbox{\boldmas #1}} %
\def\Q{\Bbb Q}                 %
\def\R{\Bbb R}                 %
\def\N{\Bbb N}                 %
\def\Z{\Bbb Z}                 %
                 %
               %
%%%%%%%%%%%%%%%%%%%%%%%%%%%%%%%%%%%%%%%%%%%%%%%%%%%%%%%%%%%%%%%%%%%%%%%

%%%%%%%%%%%%%%%%%%%%%%%%%%%%%%%%%%%%%%%%%%%%%%%%%%%%%%%%%%%%%%%%%%%%%%%
\cl{\bffg Scattered and paracompact order topologies} 
\bp
\cl{\bff Gerald Kuba}
\bp\mp
{\bff 1. Statement of main results}
\mp
Write $\,|M|\,$ for the cardinal number (the {\it size}) of a set $\,M\,$.
While an infinite set $\,S\,$ carries $\,2^{2^{|S|}}\,$ 
(mutually non-homeomorphic) topologies, $\,S\,$ 
cannot be equipped with more than $\,2^{|S|}\,$ {\it order topologies}.
This is trivial in view of the definition that $\,\tau\,$ is an 
order topology on $\,S\,$ when $\,\tau\,$ is generated by a linear 
ordering $\,\preceq\,$ of $\,S\,$. 
(A subbasis of $\,\tau\,$ is formed by all rays
$\,\{\,x\in S\;|\;x\prec y\,\}\,$ and $\,\{\,x\in S\;|\;x\succ y\,\}\,$
with $\,y\in S\,$.) In the following, a {\it linearly ordered}
space is a space whose topology is an order topology. 
Note that every linearly ordered space is hereditarily normal.
According to the title we consider only linearly ordered spaces 
which are {\it scattered}. (A Hausdorff space is scattered if and only if 
every nonempty subspace contains an isolated point.)
If $\,X\,$ is such a space then for the weight 
$\,w(X)\,$ of $\,X\,$ we always have $\,w(X)=|X|\,$.
(Because $\,w(Y)\leq|Y|\,$ trivially holds when the space $\,Y\,$ 
is linearly ordered. And it is a nice exercise to verify 
$\,|Y|\leq w(Y)\,$ for all scattered Hausdorff spaces $\,Y\,$.)
\sp
In the following we focus on scattered linearly ordered spaces 
which are either {\it metrizable} or {\it compact}.
Since size and weight coincide, a compact {\it and} metrizable 
scattered space must be countable. 
Of course, every countable compact Hausdorff space is completely
metrizable. Furthermore, a scattered metrizable space 
must be completely metrizable. 
The following theorem implies that 
every infinite set $\,S\,$ can be equipped with $\,2^{|S|}\,$
complete metrics which generate mutually non-homeomorphic 
scattered order topologies on $\,S\,$.
\mp
{\bf Theorem 1.} {\it For every cardinal $\,\kappa\geq\aleph_0\,$
there exist $\,{2^\kappa}\,$ mutually non-homeomorphic
scattered and (completely) metrizable, linearly ordered 
spaces of size $\,\kappa\,$.}
\mp
On the other hand, the following theorem implies that 
every {\it uncountable} set $\,S\,$ can be equipped with $\,2^{|S|}\,$
mutually non-homeomorphic scattered and compact order topologies.
\mp
{\bf Theorem 2.} {\it For every cardinal $\,\kappa>\aleph_0\,$
there exist $\,{2^\kappa}\,$ mutually non-homeomorphic
compact, scattered, linearly ordered spaces of size $\,\kappa\,$.}
\mp
The statement in Theorem 2 is unprovable for $\,\kappa=\aleph_0\,$
because, due to Mazurkiewicz and  Sierpi\'nski [9],
up to homeomorphism there exist precisely $\,\aleph_1\,$
compact, countable Hausdorff spaces. (They all are lineraly ordered
and scattered.) And not only $\,\aleph_1<2^{\aleph_0}\,$ 
is consistent with ZFC set theory. It is also  
consistent with ZFC  that there exist 
$\,2^{\aleph_0}\,$ cardinal mumbers $\,\kappa\,$
with $\,\aleph_0<\kappa<2^{\aleph_0}\,$.
The last statement is worth mentioning to put important 
consequences of Theorem 1 and Theorem 2 in perspective. 
Under the restriction that either $\,\kappa=\aleph_0\,$ or
$\,\kappa\geq 2^{\aleph_0}\,$
the following consequence of Theorem 1 is well-known 
and has been proved for $\,\kappa\geq 2^{\aleph_0}\,$
in the realm of pathwise connected spaces, see [2].
\mp
{\bf Corollary 1.} {\it For every cardinal $\,\kappa\geq \aleph_0\,$
there exist $\,{2^\kappa}\,$ mutually non-homeomorphic
complete metric spaces of size $\,\kappa\,$.}
\mp
Under the restriction $\,\kappa\geq 2^{\aleph_0}\,$
the following consequence of Theorem 2
has already 
been proved in the realm of connected, linearly ordered spaces, see [4].
\mp
{\bf Corollary 2.} {\it For every cardinal $\,\kappa>{\aleph_0}\,$
there exist $\,{2^\kappa}\,$ mutually non-homeomorphic
compact Hausdorff spaces of size $\,\kappa\,$.}
\mp
{\it Remark.} Not only in Theorem 2, but also in Corollary 2
the amount $\,{2^\kappa}\,$ is maximal.
(Because every compact Hausdorff space of size $\,\kappa\,$ 
is homeomorphic to a closed subspace of Hilbert cube 
$\,[0,1]^\kappa\,$.)
\bp
{\bff 2. Ultrametrics and linear orderings}
\mp
In proving Theorem 1 
we will consider scattered, {\it ultrametric} spaces only. 
(A metric $\,d\,$ on $\,X\,$ is an ultrametric if and only if 
%the non-Archimedean triangle inequality 
$\;d(x,z)\leq\max\{d(x,y),d(y,z)\}\;$
for all $\,x,y,z\in X\,$.)
By considering complete ultrametric spaces in the proof of Theorem~1 
we need not deal with linear orderings in view of the following proposition.
\mp
{\bf Proposition 1.} {\it The topology of a space is an order topology if 
it is generated by a complete ultrametric.}
\mp
One way to verify Proposition 1 is to combine two well-known results.
By [1] 6.3.2.f
the topology of any {\it strongly zero-dimensional} metrizable space 
is an order topology. And by [1] 7.3.F
a metrizable space is strongly zero-dimensional 
if and only if its topology is generated by some ultrametric.
(As a noteworthy consequence, the property {\it strongly zero-dimensional}
can be included in Theorem 1. 
In Theorem 2 the property 
{\it strongly zero-dimensional} is included automatically
because, other than for metric spaces, the properties 
{\it strongly zero-dimensional} and {\it totally disconnected}
are equivalent for compact Hausdorff spaces.)
\mp
In particular, by applying  [1] 6.3.2.f  and [1] 7.3.F,
Proposition 1 is true even without assuming completeness.
In order to keep this paper
self-contained, in the following we write down an elementary 
and direct proof of Proposition 1.
\mp
Let us call a strict linear ordering of a set $\,S\,$ a {\it DE-ordering} 
when there is a maximum and a minimum and when each point in
$\;S\setminus\{\max S\}\;$ has an immediate successor and 
each point in $\;S\setminus\{\min S\}\;$ has an immediate predecessor.
Then the corresponding order topology is discrete.
(The abbreviation {\it DE} refers to {\bf d}iscrete with  {\bf e}nd points.)
First of all we verify the following statement.
\mp
(2.1)\qquad {\it Every nonempty set can be equipped with a DE-ordering.} 
\mp
{\it Proof.} Of course, if $\,S\not=\emptyset\,$ is finite then 
every linear ordering of $\,S\,$ is a DE-ordering. 
Assume that $\,S\,$ is infinite and that $\,S\,$ is equipped with any 
well-ordering such that $\,S\,$ has a maximum. 
Then create a linearly ordered set $\,\tilde S\,$ from the 
well-ordered set $\,S\,$
by replacing each point in $\;S\setminus\{\min S,\,\max S\}\;$
with a copy of the naturally ordered set $\,\Z\,$.
Furthermore, replace $\,\min S\,$  
with a copy of the naturally ordered set $\,\N\,$
and replace $\,\max S\,$  
with a copy of the naturally ordered set $\,\Z\setminus\N\,$.
So we obtain 
a set $\,\tilde S\,$ equipped with a linear ordering 
which clearly is a DE-ordering. This is enough to verify (2.1) since $\,S\,$
and $\,\tilde S\,$ are equipollent sets.
\mp\sp
Now in order to prove Proposition 1 let $\,X\,$ be a space
whose topology is generated by an ultrametric $\,d\,$.
For every $\,n\in\N\,$ put 
\mp
\cl{$\;{\cal P}_n\,:=\,\big\{\,\{\,x\in X\;|\;d(x,a)<2^{-n}\,\}\;\,
\big|\,\;a\in X\,\big\}\,$.}
\mp
\vfill\eject
Since distinct 
ultrametrical balls with identical radii must be disjoint,
for every $\,n\in\N\,$
the family $\,{\cal P}_n\,$ is a partition of the set $\,X\,$
consisting of open-closed sets.
Since $\,B_1\subset B_2\,$ or $\,B_1\supset B_2\,$ whenever
$\,B_1,B_2\,$ are non-disjoint balls in $\,X\,$,
if $\,n>m\,$ then the partition $\,{\cal P}_n\,$ is finer
than the partition $\,{\cal P}_m\,$.
Clearly, $\;\bigcup_{n=1}^\infty{\cal P}_n\;$ is 
a basis of the space $\,X\,$. 
\sp
If the ultrametric $\,d\,$ is complete then, of course, 
$\;\bigcap_{n=1}^\infty B_n\,\not=\,\emptyset\;$
for every chain  

$\;B_1\supset B_2\supset B_3\supset\cdots\;$
of balls $\,B_n\in{\cal P}_n\,$.
\sp
For every $\,x\in X\,$ and 
for every $\,n\in\N\,$ let $\,B_n(x)\,$ denote 
the unique open-closed set $\,B\,$ with $\;x\in B\in{\cal P}_n\,$.
Then $\,B_n(x)\supset B_m(x)\,$ whenever $\,n\leq m\,$.
\mp
We define by induction 
for every $\,n\in\N\,$ a DE-ordering $\,\prec_n\,$
of the set $\,{\cal P}_n\,$ in the following way.
Firstly define a DE-ordering $\,\prec_1\,$
of the set $\,{\cal P}_1\,$ by virtue of (2.1).
\sp
If for $\,n\in\N\,$ 
a DE-ordering $\,\prec_n\,$
of the set $\,{\cal P}_n\,$ is already defined 
then define a DE-ordering $\,<_P\,$
of the set $\;\{\,A\in{\cal P}_{n+1}\;|\;A\subset P\,\}\;$
for every $\,P\in{\cal P}_n\,$ and for distinct 
$\,A,B\in{\cal P}_{n+1}\,$ put $\,A\prec_{n+1}B\,$
either when $\;A,B\subset P\in{\cal P}_n\;$ and $\,A<_P B\,$
or when $\,A\subset P\,$ and 
$\,B\subset Q\,$ and $\;P,Q\in{\cal P}_n\;$ and 
$\,P\prec_{n} Q\,$. Obviously, $\,\prec_{n+1}\,$ is a DE-ordering 
of the set $\,{\cal P}_{n+1}\,$.
\mp
Now define  $\;x\prec y\;$ for distinct $\,x,y\in X\,$  
if and only if
$\;B_n(x)\prec_n B_n(y)\;$ for some $\,n\in\N\,$. 
Then $\;x\prec y\;$ for distinct $\,x,y\in X\,$  
if and only if for some index $\,n\,$ we have 
$\;B_m(x)\prec_m B_m(y)\;$ for every index  $\,m\geq n\,$. 
Consequently, the relation $\,\prec\,$ is a strict linear ordering of the 
set $\,X\,$. We claim that 
\sp
(2.2)\quad {\it the topology of $\,X\,$ is finer than 
the order topology of $\,\prec\,$}
\sp 
and that if the ultrametric $\,d\,$ is complete then
\sp
(2.3)\quad {\it the topology of $\,X\,$ is coarser than 
the order topology of $\,\prec\,$.}
\mp 
In order to verify (2.2) it is enough to show that 
the rays 
\sp
\cl{$\;\{\,x\in X\;|\;x\prec a\,\}\;$ and
$\;\{\,x\in X\;|\;a\prec x\,\}\;$}
\sp
  are open sets in the space $\,X\,$                                          
for every $\,a\in X\,$. Let  
$\;R_a\,=\,\{\,x\in X\;|\;x\prec a\,\}\;$ for $\,a\in X\,$
and let $\,y\in R_a\,$ and choose $\,n\in\N\,$ with
$\,B_n(y)\not=B_n(a)\,$. Then $\,B_n(y)\prec_n B_n(a)\,$
and hence $\,x\prec a\,$ for every 
$\,x\in B_n(y)\,$ (since $\,B_n(x)=B_n(y)\,$ for every
$\,x\in B_n(y)\,$). Thus we have 
$\;B_n(y)\subset R_a\;$ if $\,y\in R_a\,$ and $\,B_n(y)\not=B_n(a)\,$.
In other words, every point in $\,R_a\,$ is an
interior point of $\,R_a\,$. Hence $\,R_a\,$ is open  
in the space 
$\,X\,$ for every $\,a\in X\,$. 
Similarly, $\;\{\,x\in X\;|\;a\prec x\,\}\;$ is open
in the space 
$\,X\,$ for every $\,a\in X\,$. 

\mp
In order to verify (2.3) under the completeness assumption 
it is enough to show that all sets in 
the basis $\;\bigcup_{n=1}^\infty{\cal P}_n\;$ are open with respect to the 
order topology of $\,\prec\,$. Let $\,m\in\N\,$
and $\,V\in{\cal P}_m\,$. Assume firstly that 
$\,V\,$ is neither the maximum nor the minimum 
of the DE-ordered set $\,({\cal P}_m,\prec_m)\,$
and let $\,U\,$ resp.~$\,W\,$ be the immediate predecessor resp.~successor
of $\,V\,$ in the linearly ordered set $\,({\cal P}_m,\prec_m)\,$.
For every $\,n\geq m\,$ define balls $\,U_n,W_n\,$ such that
$\,U_m=U\,$ and $\,W_m=W\,$ and that 
$\,U_{n+1}\,$ is the $\prec_{n+1}$-maximum
of $\;\{\,B\in{\cal P}_{n+1}\;|\;B\subset U_n\,\}\;$ and
that $\,W_{n+1}\,$ is the $\prec_{n+1}$-minimum 
of $\;\{\,B\in{\cal P}_{n+1}\;|\;B\subset W_n\,\}\,$.
Then for                                   
$\;\bigcap_{n=m}^\infty U_n=\{u\}\;$ and
$\;\bigcap_{n=m}^\infty W_n=\{w\}\;$ 
we obviously have 
$\;V\,=\,\{\,x\in X\;|\;u\prec x\prec w\,\}\,$.
In a similar way we obtain 
$\;V\,=\,\{\,x\in X\;|\;x\prec a\,\}\;$
for some $\,a\in X\,$ if 
$\,V\,$ is the $\prec_m$-minimum of $\,{\cal P}_m\,$
and $\;V\,=\,\{\,x\in X\;|\;a\prec x\,\}\;$
for some $\,a\in X\,$ if 
$\,V\,$ is the $\prec_m$-maximum of $\,{\cal P}_m\,$.
\mp
This concludes the proof of Proposition 1. 
\mp
{\it Remark.} If the ultrametric $\,d\,$ is not complete then (2.3) 
is not necessarily true.
For one cannot rule out the situation that 
$\,\{x\}\in{\cal P}_n\,$ for some $\,x\in X\,$ and $\,n\in\N\,$ and that 
$\,x\,$ has neither an immediate $\prec$-predecessor nor an immediate 
$\prec$-successor, whence $\,x\,$ is isolated in the space $\,X\,$ but 
not isolated with respect to the order topology. 
However, by considering special DE-orderings it is not difficult 
to accomplish (2.3) without assuming that $\,d\,$ is complete.
\bp
{\bff 3. Signature sets for metrizable spaces}
\medskip
In general topology there are two natural ways to verify that spaces 
are not homeomorphic. The first way is to apply {\it connectedness}
arguments. This way is out of the question  
when we deal with totally disconnected spaces. 
The second way is to use {\it Cantor derivatives}.
Let $\,\Omega\,$ denote the class of all ordinals, whence
$\;\{0\}\cup\N\,\subset\,\Omega\,$.
If $\,X\,$ is a Hausdorff space and $\,\xi\in\Omega\,$
and $\,A\subset X\,$ then 
$\,A^{(\xi)}\,$ is the $\xi$-{\it th derivative} of the point set $\,A\,$.
($\,A^{(0)}=A\,$ and $\,A'\,$ is the set of all limit points of $\,A\,$ and if
$\;\alpha\in\Omega\;$ then $\,A^{(\alpha+1)}={A^{(\alpha)}}'\,$ 
and if $\,\lambda>0\,$ is a limit ordinal then
$\;A^{(\lambda)}\,:=\,\bigcap\,\{\,A^{(\alpha)}\;|\;
0<\alpha<\lambda\,\}\,$.) 
Naturally, if $\,0\not=\xi\in\Omega\;$ then $\,A^{(\xi)}\,$ is closed. 
Furthermore, $\,X\,$ is scattered if and only if 
$\,X^{(\alpha)}=\emptyset\,$ for some $\,\alpha\in \Omega\,$.
\medskip\sp
{\bf Lemma 1.} {\it Let $\,Z\,$ be a Hausdorff space with 
$\,Z''=\emptyset\,$
and let $\,H\,$ be a Hausdorff space. Then for the product space
$\;Z\times H\;$ we have 
$\;(Z\times H)^{(\xi+1)}\;=\;
Z'\times H^{(\xi)}\;\cup\; Z\times H^{(\xi+1)}\;$ 
for every ordinal number $\,\xi\,$.}
\sp
{\it Proof.} Let $\,{\bf E}[\xi]\,$ denote the equation 
$\;(Z\times H)^{(\xi+1)}\;=\;
Z'\times H^{(\xi)}\;\cup\; Z\times H^{(\xi+1)}\,$. 
Similarly as in calculus,
the product formula $\;(X\times Y)'\,=\,X'\!\times\! Y\,\cup\,X\!\times\! Y'\;$
is true for arbitrary Hausdorff spaces $\,X,Y\,$. 
(For $\,(x,y)\,$ is isolated in $\,X\times Y\,$ if and only if 
$\,x\,$ is isolated in $\,X\,$ and $\,y\,$ is isolated in $\,Y\,$.)
This has two consequences.
Firstly, $\,{\bf E}[\xi]\,$ is true for $\,\xi=0\,$.
Secondly if $\,{\bf E}[\xi]\,$ holds for $\,\xi=\alpha\,$ then
$\,{\bf E}[\xi]\,$ holds for $\,\xi=\alpha+1\,$.
The set $\;Z'\times H^{(\xi)}\;\cup\; Z\times H^{(\xi+1)}\;$
is the union of the disjoint sets 
$\;Z'\times H^{(\xi)}\;$ and $\;(Z\setminus Z')\times H^{(\xi+1)}\;$
since $\;H^{(\xi+1)}\subset H^{(\xi)}\,$.
Therefore, if $\,\lambda>0\,$ is a limit ordinal and 
$\,{\bf E}[\xi]\,$ is true for every $\,\xi<\lambda\,$ 
then $\;(Z\times H)^{(\lambda)}\;=\;
\bigcap_{0<\xi<\lambda}Z'\times H^{(\xi)}\;\cup\;
\bigcap_{0<\xi<\lambda}(Z\setminus Z')\times H^{(\xi+1)}\;
=\;Z'\times H^{(\lambda)}\;\cup\; (Z\setminus Z')\times H^{(\lambda)}\;
=\;Z\times H^{(\lambda)}\;$
and hence 
$\;(Z\times H)^{(\lambda+1)}\;=\;((Z\times H)^{(\lambda)})'\;=\;
(Z\times H^{(\lambda)})'\;=\;
Z'\times H^{(\lambda)}\;\cup\; Z\times H^{(\lambda+1)}\;$
by applying the product formula,
whence $\,{\bf E}[\xi]\,$ is true for $\,\xi=\lambda\,$, {\it q.e.d.}
\mp\sp
If $\,X\,$ is a Hausdorff space and $\,\kappa\,$ an uncountable cardinal 
then let $\,C_\kappa(X)\,$ denote the set of all points $\,x\in X\,$ such that 
$\;|U|\geq\kappa\;$ for every neighborhood of $\,x\,$ and 
$\;|U|=\kappa\;$ for some neighborhood of $\,x\,$.
(One may call the members of $\,C_\kappa(X)\,$ the 
$\kappa$-{\it condensation points} of $\,X\,$.)
Define a signature set $\,\Sigma[\kappa;X]\,$ by
\mp
\centerline{$\Sigma[\kappa;X]\;:=\;
\big\{\,\alpha\,\in\,\Omega\setminus\{0\}\;\,\big|\,\;
(X^{(\alpha)}\setminus X^{(\alpha+1)})\cap C_\kappa(X)  
\,\not=\,\emptyset\,\big\}\;.$}
\mp
Of course, the class $\,\Sigma[\kappa, X]\,$ is always a set
and two spaces $\,X_1\,$ and $\,X_2\,$ cannot be homeomorphic if 
$\;\Sigma[\kappa;X_1]\not=\Sigma[\kappa;X_2]\;$
for some cardinal $\,\kappa>\aleph_0\,$. 
\mp
For the proof of Theorem 1 we also need another signature set.
For any Hausdorff space 
$\,X\,$ let $\,\Gamma(X)\,$ be the set of all points $\,x\in X\,$
such that no neighborhood of $\,x\,$ is compact
(or, equivalently, the closure of an open neighborhood
of $\,x\,$ is never  compact).
Furthermore, define 

\centerline{$\Sigma(X)\;\,:=\,\;\big\{\,k\in{\Bbb N}\;\,\big|\,\;
\big(X^{(k)}\setminus X^{(k+1)}\big)\cap \Gamma(X)\,\not=
\,\emptyset\,\big\}$}
\mp
Note that we regard $\,\N\,$ to be defined in the classical way, 
i.e.~$\,\N\,$ does not contain $0$.
Clearly, two spaces $\,X_1\,$ and $\,X_2\,$ cannot be homeomorphic if 
$\;\Sigma(X_1)\not=\Sigma(X_2)\,$.

\vfill\eject
\bp
{\bff 4. Countable Polish spaces} 
\mp
Let $\,c=2^{\aleph_0}\,$ denote the cardinality of the continuum.
The size of a perfect
completely metrizable space cannot be smaller than $\,c\,$
(cf.~[1] 4.5.5). Consequently, all completely metrizable 
spaces of size smaller than $\,c\,$ are scattered.
Furthermore, a Polish space 
is scattered if and only if it is countable.
In view of this fact and by virtue of Proposition 1
the special case $\,\kappa=\aleph_0\,$ in Theorem 1 
is settled by the following observation.
\sp
(4.1) $\,$ {\it There exist $\,c\,$ mutually non-homeomorphic 
closed and countable subspaces of the Polish space $\,\R\setminus\Q\,$.}
\sp
Notice that Proposition 1 can be applied because 
the topology of every space provided by (4.1) 
can be generated by a complete ultrametric. 
This is a consequence of the observation that 
the topology of the Baire space $\,\R\setminus\Q\,$
is generated by a very natural complete ultrametric.
(Declare $\,2^{-n}\,$ as the distance of distinct irrationals 
$\,a\,$ and $\,b\,$ when $\,n\,$ is the smallest 
index of distinct quotients in the continued fractions 
of $\,a\,$ and $\,b\,$.) 
\mp
{\it Remark.} 
While for both $\;X=\R\;$ and $\;X\,=\,\R\setminus\Q\;$
the order topology on $\,X\,$ generated by the natural ordering of 
the elements of $\,X\,$ coincides with the Euclidean topology on $\,X\,$,
such a coincidence is not true for subsets of $\,X\,$.
Even worse, while such a coincidence holds
for {\it closed} subsets of the connected space $\,\R\,$, 
it does not necessarily hold for closed 
subsets of the totally disconnected Baire space  $\,\R\!\setminus\!\Q\,$.
(For example, $\;\{-e\}\cup\{\,e^{-n}\;|\;n\in\N\,\}\;$ is an infinite,
closed, discrete subspace of $\,\R\!\setminus\!\Q\,$ whose order topology 
is compact. See also Proposition 2 below.)
Therefore, Proposition 1 is essential for deriving the 
case $\,\kappa=\aleph_0\,$ in Theorem 1 from (4.1).
\mp
The statement
(4.1) equals the classic solution of the enumeration problem 
concerning {\it countable Polish spaces} due to 
Mazurkiewicz and  Sierpi\'nski (see [2] Lemma 4.10).
In order to keep this paper
self-contained we will prove (4.1) and hence 
Theorem 1 for $\,\kappa=\aleph_0\,$
in this section. And by proving a bit more we will gain a deeper insight
in connection with the following observation.
\sp
(4.2) $\,$ {\it It is unprovable that $\,{\Bbb R}\,$ 
contains $\,c\,$ mutually non-homeomorphic closed 
and countable subspaces.}
\sp
The observation (4.2) is true because by [9] 
there exist precisely $\,\aleph_1\,$ compact and countable subspaces
of $\,\R\,$ up to homeomorphism and therefore (see [5] Theorem 8.1)
the real line $\,\R\,$ has precisely $\,\aleph_1\,$ 
closed and countable subspaces up to homeomorphism.
Motivated by comparing (4.1) and (4.2) 
we are now going to prove the following proposition which 
implies (4.1) and hence settles the case $\,\kappa=\aleph_0\,$ in Theorem~1.
Throughout this section, if $\,X\subset{\Bbb R}\,$ then
$\,\overline{X}\,$ denotes the space which equals 
the closure of $\,X\,$ in the real line $\,{\Bbb R}\,$
equipped with the Euclidean topology, and  
$\,\hat X\,$ denotes the space which equals the set $\,X\,$ 
equipped with the order 
topology generated by the natural ordering of the reals in $\,X\,$.
\medskip
{\bf Proposition 2.} {\it There exists a family $\;{\cal Y}\;$
of mutually non-homeomorphic countable subspaces of $\,{\Bbb R}\,$ 
such that $\;|{\cal Y}|=c\;$ and 
$\;\overline{Y}\setminus {\Bbb Q}\,=\,Y\;$ for every $\,Y\in{\cal Y}\,$
and all spaces in the family 
$\;\{\,\overline{Y}\;|\;Y\!\in\!{\cal Y}\,\}\cup
\{\,\hat Y\;|\;Y\!\in\!{\cal Y}\,\}\;$
are homeomorphic.}
\mp
{\it Proof.} In the following put $\,\N^*=\N\setminus\{1\}\,$,
whence $\;\N^*\,=\,\{\,k\in\Z\;|\;k\geq 2\,\}\,$.
For $\,n\in{\Bbb N}\,$ 
let $\,K_n\,$
be a compact, countable, well-ordered subset of 
$\;[\pi+2n,\pi+2n+1]\setminus \Q\;$
such that $\;K_n^{(n)}=\{\max K_n\}=\{\pi+2n+1\}\,$.
(The set $\,K_n\,$ may be defined as an appropriate 
order-isomorphic copy of the canonically ordered set of all ordinals
$\,\alpha\leq\omega^n\,$. It is also straightforward to construct 
$\;K_1,\,K_2,\,K_3,\,...\;$ recursively 
without using ordinal numbers.)
\eject
\mp
In order to prove Proposition 2, let $\;r_0,r_1,r_2,...\;$ 
be a strictly decreasing 
sequence of {\it rational} numbers in $\,[\pi,\pi+1]\,$
with $\;\inf\,\{\;r_m\;|\;m\in {\Bbb N}\,\}\,=\,\pi\,$.
For every $\,m\in{\Bbb N}\,$ let $\,\xi_m(1),\xi_m(2),\xi_m(3),...\,$
be a strictly decreasing 
sequence of {\it irrational} numbers in $\,[r_{m},r_{m-1}]\,$
with $\;\inf\,\{\;\xi_m(k)\;|\;k\in {\Bbb N}\,\}\,=\,r_{m}\,$.
For every $\,n\in{\Bbb N}\,$ define a {\it discrete} subset $\,E_n\,$ of 
$\;]\pi+2n+1,\pi+2n+2[\,\setminus {\Bbb Q}\;$ via
$\;E_n\,:=\,
\big\{\,2n+1+\xi_m(k)\;\,\big|\,\;m,k\in{\Bbb N}\,\big\}\,$.
\mp
For every infinite subset $\,S\,$ of $\,\N^*\,$ put 
\medskip
\centerline{$\;\;G_S\;:=\;\bigcup\limits_{n\in S}(K_n\cup E_n)\;$.}
\medskip  
For every $\,n\in{\Bbb N}\,$ we have
$\;\overline{E_n}\setminus E_n\,
=\,\{\,2n+1+r_m\;|\;m\in{\Bbb N}\,\}\cup\{\pi+2n+1\}\;$
and hence the only irrational limit point of $\,E_n\,$ is 
$\;\pi+2n+1\,=\,\max  K_n\,$. 
So for every infinite subset $\,S\,$ of $\,{\Bbb N}^*\,$
we have $\;G_S\,=\,\overline{G_S}\setminus{\Bbb Q}\;$
since $\;\overline{G_S}\;=\;G_S\cup\bigcup_{n\in S}\overline{E_n}\,$.
We claim that the family 
\medskip
\centerline{$\;{\cal Y}\;:=
\;\{\,G_S\;\,|\,\;S\subset{\Bbb N}\;\land\;|S|=\aleph_0\,\}\;$}
\medskip
is as desired. Obviously, we always have 
$\;\Gamma(G_S)\,=\,\{\,\max  K_n\;|\;n\in S\,\}\,$.
Therefore, since the point $\,\max K_n\,$ lies 
in $\;(K_n\setminus\{\max K_n\})^{(m)}\;$ but not in 
$\;(K_n\setminus\{\max K_n\})^{(m+1)}\;$ if and only if $\,m=n\,$,
we must always have
$\;\Sigma(G_S)=S\,$. Hence two spaces $\,G_S\,$ and $\,G_T\,$
are never homeomorphic for distinct infinite sets $\;S,T\subset{\Bbb N}^*\,$.
\mp
We conclude the proof by showing that 
for some space $\,W\,$ 
both spaces $\,\overline{Y}\,$ and $\,\hat Y\,$ are 
homeomorphic to $\,W\,$ for each $\,Y\in{\cal Y}\,$.
Let $\,S\,$ be an arbitrary infinite subset of $\,\N^*\,$.
The order-type of the naturally ordered set  
$\;\{\,-x\;|\;x\in\overline E_n\,\}\;$ is $\,\omega^2+1\,$
while the order type of $\,K_n\,$ is $\,\omega^n+1\,$.
Hence by a {\it zipper} argument (in a Hilbert's hotel kind of way)
it is plain to find a homeomorphism from the compact space
$\;K_n\cup\overline{E_n}\;$ onto the compact space $\,K_n\,$
for every $\,n\geq 2\,$. Consequently, since each compact building block 
$\; K_n\cup \overline{E_n}\;\,(n\in S)\;$
is open in the space
$\;\overline{G_S}\,=\,\bigcup_{n\in S}( K_n\cup \overline{E_n})\,$,
the space $\,\overline{G_S}\,$ 
is homeomorphic to $\;V_S\,:=\,\bigcup_{n\in S} K_n\,$.
The naturally ordered set 
$\,V_S\,$ is well-ordered without a maximum
and $\;V_S^{(k)}\not=\emptyset\;$ for every $\,k\in{\Bbb N}\,$
but $\;\bigcap_{k=1}^\infty V_S^{(k)}=\emptyset\,$.
Hence $\,V_S\,$ is order-isomorphic to the well-ordered set $\,W\,$
of all ordinal numbers smaller than $\,\omega^\omega\,$.
Therefore, since for every closed subset $\,A\,$ of 
$\,{\Bbb R}\,$ the linearly ordered space 
$\,\hat A\,$ is identical with the subspace $\,A\,$ of the real line,
the space $\,\overline{G_S}\,$ is homeomorphic 
to the  space $\,W\,$ equipped with the canonical order topology.
Finally it is evident that 
the linearly ordered space $\,\hat G_S\,$
is homeomorphic to the 
subspace $\,\overline{G_S}\,$ of $\,\R\,$ and hence to 
the space $\,W\,$, {\it q.e.d.}
\bp\sp
{\bff 5. Ultrametric hedgehogs}
\mp
Let us call a space 
{\it completely ultrametrizable} if and only if 
its topology is generated by some complete ultrametric.
(It is worth mentioning that a space must be 
completely ultrametrizable
if its topology is generated by some ultrametric and
also by some complete metric, see [7] Lemma 3.
Therefore, since every scattered metrizable space is completely metrizable,
in view of Proposition 1 it would be enough 
to consider ultrametrics instead of complete ultrametric.
But the ultrametrics we will consider in the proof of Theorem 1 
are very natural 
and immediately recognized as complete metrics.)
\sp
Clearly, if $\,X_1\,$ and $\,X_2\,$ are 
completely ultrametrizable spaces then the product space $\;X_1\times X_2\;$
and the topological sum of the two spaces
are completely ultrametrizable as well. 
(Note that if $\,d_i\,$ is an ultrametric on $\,X_i\,$
then the standard {\it maximum metric} on $\;X_1\times X_2\;$ 
with respect to $\,d_1,d_2\,$ is an ultrametric.)
\vfill\eject
\mp
For an index set $\,I\not=\emptyset\,$ let $\,X_i\,$ be an infinite 
metric space and $\,b_i\,$ a point in $\,X_i\,$ for every 
$\,i\in I\,$. The metric of $\,X_i\,$ is denoted by $\,d_i\,$
and we assume that $\,d_i\,$ is an ultrametric.
  Fix $\,o\,$ such that 
$\;o\,\not\in\,\bigcup_{i\in I}I\times X_i\;$ and put 
\mp
\cl{$\;{\cal H}[(X_i,b_i)_{i\in I}]\;=\;
\{o\}\;\cup\,\bigcup\limits_{i\in I}\{i\}\times(X_i\setminus\{b_i\})\;$.}
\mp
Define a function $\,d\,$ from $\,{\cal H}[(X_i,b_i)_{i\in I}]^2\,$
into $\,\R\,$ as follows: 
\sp
Put $\;d(o,o)=0\,,$ and for each $\,i\in I\,$ put
\sp
$\;d\big((i,x),o\big)=d\big(o,(i,x)\big)=d_i(b_i,x)\;$ 
whenever $\,b_i\not=x\in X_i\,$,  
and for each $\,i\in I\,$ put
\sp
$\;d\big((i,x),(i,y)\big)=d_i(x,y)\;$ 
whenever $\;x,y\,\in\,X_i\setminus\{b_i\}\,$, and 
for distinct $\,i,j\in I\,$ put
\sp
$\;d\big((i,x),(j,y)\big)\,=\,\max\{d_i(b_i,x),d_j(b_j,y)\}\;$
whenever $\,b_i\not=x\in X_i\,$ and $\,b_j\not=y\in X_j\,$.
\mp
After having thoroughly checked that in all possible scenarios 
any triangle is isosceles
and any non-equilateral triangle has one side shorter
than the other sides, we see that 
$\,d\,$ is an ultrametric on the set $\,{\cal H}[(X_i,b_i)_{i\in I}]\,$.
Trivially, the subspace 
$\;\{o\}\cup(X_i\setminus\{b_i\})\;$ of
$\,{\cal H}[(X_i,b_i)_{i\in I}]\,$ is an isometric copy of $\,X_i\,$
for every $\,i\in I\,$. 
(So if the sets $\,X_i\,$ are mutually disjoint
then every basic space $\,X_i\,$ can be identified with 
the space $\,\{o\}\cup(X_i\setminus\{b_i\})\,$ 
in order to achieve $\;{\cal H}[(X_i,b_i)_{i\in I}]\,=\,\bigcup_{i\in I}X_i\;$
and $\;X_i\cap X_j\,=\,\{o\}\;$ for all distinct $\,i,j\in I\,$.)
By analogy to the classical hedgehog of spininess $\,|I|\,$ (see [1] or [10])
we call $\;{\cal H}[(X_i,b_i)_{i\in I}]\;$
an {\it ultrametric hedgehog}
and we call the basic spaces $\;X_i\;(i\in I)\;$ the {\it spines}
and the point $\,o\,$ the {\it body} 
of the hedgehog.
\mp
It is evident that the hedgehog-metric $\,d\,$ is complete
if and only if all spine-metrics $\,d_i\,$ are complete.
The weight of $\,{\cal H}[(X_i,b_i)_{i\in I}]\,$
is not smaller than $\,\max\{|I|,\aleph_0\}\,$ because each spine contains an 
open set disjoint from $\,\{o\}\,$.
If $\,w(X_i)=\kappa\geq\aleph_0\,$ for every $\,i\in I\,$
then $\;w({\cal H}[(X_i,b_i)_{i\in I}])=\,\max\{|I|,\kappa\}\;$
because each spine contains a dense set of size $\,\kappa\,$.
\bp
{\bff 6. Proof of Theorem 1}
\mp\sp
{\bf Proposition 3.} {\it For every cardinal $\,\kappa>\aleph_0\,$ there 
exists a complete ultrametric space $\,Y_\kappa\,$ of size $\,\kappa\,$
such that $\,Y_\kappa'=C_\kappa(Y_\kappa)\,$ and 
$\,|C_\kappa(Y_\kappa)|=1\,$ and hence
$\,Y_\kappa\setminus C(Y_\kappa)\,$ is a discrete subspace of $\,Y_\kappa\,$.}
\mp
{\it Proof.} Consider the set 
$\;Z\,:=\,\{0\}\cup\{\,2^{-n}\;|\;n\in\N\,\}\;$ 
and declare $\;\max\{x,y\}\;$ as the distance between distinct
points $\,x,y\in Z\,$. It is evident that this distance 
defines a complete ultrametric on $\,Z\,$ which generates the Euclidean 
topology restricted to $\,Z\,$. Thus $\;Z'\,=\,\{0\}\,$.
Let $\,I\,$ be an index set of size $\,\kappa>\aleph_0\,$
and for every $\,i\in I\,$ let $\,A_i\,$ be identical 
with the complete ultrametric space $\,Z\,$.
Put $\;Y_\kappa\,:=\,{\cal H}[(A_i,0)_{i\in I}]\;$
and let $\,o\,$ be the body of this ultrametric hedgehog.
We observe that $\,Y_\kappa'=C_\kappa(Y_\kappa)=\{o\}\,$, {\it q.e.d.}
\mp\sp
{\bf Proposition 4.} {\it For every ordinal $\,\alpha>0\,$
there exists a 
complete ultrametric space $\,Z_\alpha\,$ 
of size $\,\max\{\aleph_0,|\alpha|\}\,$ 
such that $\,Z_\alpha^{(\alpha)}\,$ is a singleton.
(In particular, $\,Z_\alpha\,$ is scattered.)}
\mp
{\it Proof}. Of course,  we construct the desired spaces
by transfinite induction. 
Put $\,Z_1=Z\,$ where $\;Z\,=\,\{0\}\cup\{\,2^{-n}\;|\;n\i\N\,\}\;$ 
is the complete ultrametric space as in the proof of Proposition 3.
Suppose firstly that for an ordinal $\,\alpha>0\,$ a complete ultrametric
space $\,Z_\alpha\,$ with $\;Z_\alpha^{(\alpha)}=\{z_\alpha\}\;$ 
and $\;|Z_\alpha|=\max\{\aleph_0,|\alpha|\}\;$ 
is already defined. (This is true for $\,\alpha=1\,$.)
Then we define $\,Z_{\alpha+1}\,$ as the product space $\;Z\times Z_\alpha\;$
and consider it equipped with the maximum metric
which is a complete ultrametric. Clearly, 
$\;|Z_{\alpha+1}|=\max\{\aleph_0,|\alpha|\}=\max\{\aleph_0,|\alpha+1|\}\,$.
By virtue of Lemma 1 we have $\;Z_{\alpha+1}^{(\alpha+1)}=\{(0,z_\alpha)\}\,$.
\sp
Secondly, let $\,\lambda>0\,$ be a limit ordinal
and assume that 
for every ordinal $\,\alpha\,$ with $\,1\leq \alpha<\lambda\,$
a complete ultrametric
space $\,Z_\alpha\,$ with $\;Z_\alpha^{(\alpha)}=\{z_\alpha\}\;$ 
and $\;|Z_\alpha|=\max\{\aleph_0,|\alpha|\}\;$ 
is already defined.
Then let $\,Z_\lambda\,$ be the complete ultrametric 
hedgehog $\;{\cal H}[(Z_\alpha,z_\alpha)_{1\leq \alpha<\lambda}]\,$.
Obviously, the body $\,o\,$ of the hedgehog $\,Z_\lambda\,$
lies in $\,Z_\lambda^{(\beta)}\,$ whenever $\,1\leq \beta<\lambda\,$.
Furthermore, $\;Z_\alpha^{(\beta)}=\emptyset\;$
whenever $\,1\leq \alpha<\beta<\lambda\,$.
Consequently, 
$\;Z_\lambda^{(\lambda)}\,=\,\bigcap_{1\leq \beta<\lambda}Z_\lambda^{(\beta)}
\,=\,\{o\}\,.$ 
Of course, $\;|Z_\lambda|=|\lambda|=\max\{\aleph_0,|\lambda|\}\,$.
This concludes the proof of Proposition 4.
\bp
Now we are ready to prove Theorem 1.
The special case $\,\kappa=\aleph_0\,$ is settled by Propositions 1 and 2.
Therefore we assume $\,\kappa>\aleph_0\,$. Let
$\;K\,:=\,\{\,\alpha\in\Omega\;|\;\omega\leq \alpha <\kappa\,\}\,.$
Then $\;|K|=\kappa\;$ and hence $\,K\,$ has precisely 
$\,2^\kappa\,$ nonempty subsets.
Let $\,Y_\kappa\,$ be an ultrametric space 
with $\,|C_\kappa(Y_\kappa)|=1\,$
as in the proof of Proposition 3
and for each $\,\alpha\in K\,$ let $\,Z_\alpha\,$ be an 
ultrametric space as in Proposition 4 with 
$\,Z_\alpha^{(\alpha)}=\{z_\alpha\}\,$.
Since $\,0\not\in K\,$ we may unambiguously 
write $\,Z_0\,$ for the space $\,Y_\kappa\,$ and $\,z_0\,$
for the unique condensation point in $\,Y_\kappa\,$.
For every $\,\alpha\in K\,$ put 
$\;H_\alpha\,:=\,{\cal H}[(Z_\beta,z_\beta)_{\beta\in\{0,\alpha\}}]\,$.
So  $\,H_\alpha\,$ is a complete ultrametric 
hedgehog with the two spines $\,Y_\kappa\,$ and $\,Z_\alpha\,$. 
Clearly, $\,|H_\alpha|=\kappa\,$ and 
if $\,o\,$ is the body of the hedgehog $\,H_\alpha\,$ then 
$\;H_\alpha^{(\alpha)}=C_\kappa(H_\alpha)=\{o\}\,$.
(Note that $\;|Z_\alpha|<\kappa\;$ for every $\,\alpha\in K\,$.)
In view of the proofs of Propositions 3 and 4
we may assume that the distance 
between points in $\,H_\alpha\,$ is always smaller than $1$.
(Alternatively, if $\,(X,d)\,$ is a complete ultrametric space then 
it is evident 
that $\;\min\{{8\over 9},d(\cdot,\cdot)\}\;$ is a complete
ultrametric which generates the topology of $\,(X,d)\,$.)
\sp
Finally, for each nonempty $\,L\subset K\,$ define $\,{\cal S}[L]\,$ 
as the {\it standard metrical sum} of the 
complete ultrametric spaces
$\;H_\alpha\;(\alpha\in L)\,$. 
(The distance between any two points in a common summand 
remains unchanged and the distance between two points in
distinct summands always equals $1$.)
Clearly, $\,{\cal S}[L]\,$ is a scattered and 
complete ultrametric space of size $\,\kappa\,$.
The $\,2^\kappa\,$ spaces
$\;{\cal S}[L]\;(\emptyset\not=L\subset K)\;$
are mutually non-homeomorphic because we obviously have  
$\;\Sigma[\kappa;{\cal S}[L]]=L\;$ 
whenever $\;\emptyset\not=L\subset K\,$.
\bp
{\bff 7. Signature sets for non-metrizable spaces}
\medskip
In order to prove Theorem 2 we also work with Cantor derivatives.
Again let $\,\Omega\,$ be the canonically well-ordered 
class of all ordinal numbers (with $\,{\Bbb N}\cup\{0\}\,\subset\,\Omega\,$).
Furthermore let $\,{\cal L}\,$ denote the class of all {\it limit ordinals}
(with $\,0\in {\cal L}\,$).
Note that any nonempty sub{\it set} of the {\it class} $\,\Omega\,$
has a well-defined {\it supremum.} Put
$\;[\alpha,\beta]\,:=\,\{\,\xi\in\Omega\;|\;\alpha\leq\xi\leq\beta\,\}\;$ 
and $\;[\alpha,\beta[\,:=\,[\alpha,\beta]\setminus\{\beta\}\;$
and $\;]\alpha,\beta[\,:=\,[\alpha,\beta]\setminus\{\alpha,\beta\}\;$
whenever $\,\alpha,\beta\in \Omega\,$. 
If we speak of the {\it space} 
$\,[\alpha,\beta]\,$ resp.~$\,[\alpha,\beta[\,$
resp.~$\,]\alpha,\beta[\,$ then 
we refer to the order topology of the canonical well-ordering.
Note that all spaces $\,[\alpha,\beta]\,$ are compact
and scattered and hereditarily normal. 
\medskip
If $\,\alpha\in\Omega\,$ then put $\;|\alpha|\,:=|[0,\alpha[|\,$. 
(This definition is a tautology if ordinal numbers are defined
in the standard way.) Put
$\;o(\kappa)\,:=\,\min\big\{\,\gamma\in\Omega\;\,\big|\,\;
|\gamma|=\kappa\,\big\}\;$
for every cardinal number $\,\kappa\,$.
(If cardinal numbers are defined as initial ordinal 
numbers then, of course, $\,o(\kappa)=\kappa\,$ 
for every $\,\kappa\,$.)
For any cardinal number $\,\kappa\,$ 
let (as usual) $\,\kappa^+\,$ denote the smallest
cardinal number greater than $\,\kappa\,$. 
(For example,
$\,\aleph_1=\aleph_0^+\,$.)
For cardinals $\,\kappa\,$ and ordinals $\,\alpha\,$ 
we have $\;|\alpha|=\kappa\;$ if and only if 
$\;o(\kappa)\leq\alpha<o(\kappa^+)\,$. 
In particular, 
$\,\omega:=o({\aleph_0})\,$ is the smallest infinite ordinal
and $\,\omega_1:=o({\aleph_1})\,$ is the smallest uncountable ordinal.
\medskip
For $\,\xi\in\Omega\,$ we write (as usual) $\,\omega^\xi\,$
for the {\it ordinal power} with basis $\,\omega\,$ and exponent $\,\xi\,$.
So all spaces $\,[0,\omega^\xi]\,$ are compact and
for $\,\xi>0\,$ we have $\;|[0,\omega^\xi]|=\max\{\aleph_0,|\xi|\}\,$.
In particular, $\;|[0,\omega^\xi]|=|\xi|\;$ for every ordinal 
$\,\xi\geq \omega\,$.
\medskip\smallskip
As above, if $\,X\,$ is a Hausdorff space and $\,\xi\in\Omega\,$
and $\,A\subset X\,$ then 
$\,A^{(\xi)}\,$ is the $\xi$-{\it th derivative} of the point set $\,A\,$.
Clearly, $\;A^{(\alpha)}\supset A^{(\beta)}\;$ whenever
$\;0<\alpha\leq\beta\;$ and for $\;A\subset B\subset X\;$
we have $\;A^{(\alpha)}\subset B^{(\alpha)}\;$ for every 
$\,\alpha\in\Omega\,$.
The following lemma is evident.
(Historically, Cantor's definition of 
the ordinal powers of $\,\omega\,$ is designed precisely 
so that the following is true.)
\medskip
{\bf Lemma 2.} {\it Let $\,0\not=\xi\in \Omega\,$.
With respect to the compact space $\,[0,\omega^\xi]\,$,
for every ordinal $\,\alpha>0\,$ 
the point sets $\,[0,\omega^\xi[^{(\alpha)}\,$ and
$\,[0,\omega^\xi]^{(\alpha)}\,$ coincide and they 
contain the point $\,\omega^\xi\,$ if and only if 
$\;\alpha\leq \xi\,$. And 
$\;[0,\omega^\xi]^{(\xi)}=[0,\omega^\xi[^{(\xi)}=\{\omega^\xi\}\,$.}
\bigskip
In the following we distinguish between 
{\it regular} and {\it singular} cardinal numbers. 
Singular cardinals are the cardinals which are not regular. 
A cardinal number $\,\kappa\,$ is {\it regular}
if and only if $\;\sup A\,<\,o(\kappa)\;$
for every subset $\,A\,$ of 
$\;\{\,\alpha\in \Omega\;|\;\alpha<o(\kappa)\,\}\;$
with $\,|A|<\kappa\,$. (Note that any sub{\it set}
of the class $\,\Omega\,$ has a well-defined supremum in $\,\Omega\,$.)
Topologically speaking, 
{\it an infinite cardinal number $\,\kappa\,$ is 
regular if and only if in the compact linearly ordered space 
$\;[0,o(\kappa)]\;$
the first derivative of a point set $\,A\,$ with $\,|A|<\kappa\,$ 
does not contain the point $\,o(\kappa)\,$.}
(In the following it is essential 
that $\,\kappa^+\,$ is regular for every $\,\kappa\,$.) 
\medskip
If $\,(X,\preceq)\,$ is a linearly ordered set then let 
$\;[a,b]\,=\,\{\,x\in X\;|\;a\preceq x\preceq b\,\}\;$
and $\;]a,b[\;=\,\{\,x\in X\;|\;a\prec x\prec b\,\}\;$
and $\;[a,b[\;=\,[a,b]\setminus\{b\}\;$
for $\,a,b\in X\,$. Furthermore let $\,\preceq^*\,$ denote
the {\it backwards linear ordering} defined by 
$\,x\prec^* y\,$ if and only if $\,y\prec x\,$.
In the usual sloppy way, if $\,X\,$ is a set of ordinal numbers,
then $\,X^*\,$ 
is the set $\,X\,$ equipped with the backwards linear ordering of 
the canonical well-ordering of $\,\Omega\,$.
If $\,(X,\preceq)\,$ and $\,(Y,\leq)\,$
are two linearly ordered sets then the 
{\it lexicographic ordering} of any nonempty subset 
$\,Z\,$ of $\,X\times Y\,$ is defined so that 
$\,(x_1,y_1)\,$ is smaller 
than $\,(x_2,y_2)\,$ when either 
$\;x_1\prec x_2\;$ or when $\;x_1=x_2\;$
and $\,y_1<y_2\,$. 
\medskip
In the following let $\,X\,$ be a {\it scattered} 
Hausdorff space and $\,\kappa>\aleph_0\,$ be
a regular cardinal number. As in Section 2 let
$\,C_\kappa(X)\,$ denote
the set of all $\kappa$-condensation points.
(So a point $\,x\in X\,$ lies in $\,C_\kappa(X)\,$ 
if and only if  $\,|U|\geq \kappa\,$
for every neighborhood $\,U\,$ of $\,x\,$
and $\,|U|=\kappa\,$ for some neighborhood $\,U\,$ of $\,x\,$.)
For example,  $\;C_{\kappa}([0,o(\kappa)])\,=\,
\{o(\kappa)\}\;$ for every regular cardinal $\,\kappa> \aleph_0\,$.
For $\,x\in X\,$ let 
$\,\Omega_\kappa(x)\,$ denote the class of all ordinals $\,\alpha\,$
such that $\;x\in A^{(\alpha)}\;$ for some point set 
$\,A\subset X\,$ with $\;|A|<\kappa\,$. The class $\,\Omega_\kappa(x)\,$
is never empty since, trivially, 
$\;0\in \Omega_\kappa(x)\;$ for every $\,x\in X\,$.
Since $\,X\,$ is scattered, the class $\,\Omega_\kappa(x)\,$ is a
nonempty set for each $\,x\in X\,$. Moreover, 
$\;\Omega_\kappa(x)\subset [0,o(\kappa)[\;$
because $\;A^{(o(\kappa))}=\emptyset\;$ whenever $\,A\subset X\,$
and $\,|A|<\kappa\,$.
\medskip
So we may define a signature set with respect to 
the scattered space $\,X\,$ and the regular cardinal $\,\kappa\,$ by
\medskip
\centerline{$\;{\Psi}[X,\kappa]\;:=\;
\{\,\sup \Omega_\kappa(x)\;|\;x\in C_\kappa(X)\,\}\,$.}
\medskip\smallskip
Clearly, two scattered spaces $\,X_1,X_2\,$ cannot be homeomorphic 
if $\;{\Psi}[X_1,\kappa]\not={\Psi}[X_2,\kappa]\;$
for some regular cardinal $\,\kappa\,$.
For each regular cardinal $\,\kappa> \aleph_0\,$
we have $\;{\Psi}[[0,o(\kappa)],\kappa]=\{0\}\,$.
More generally, in view of the following lemma,
$\;{\Psi}[[0,\alpha],\kappa]\,\subset\,
\{0,o(\kappa)\}\;$ for every $\,\alpha\in\Omega\,$.
\bigskip
{\bff 8. Proof of Theorem 2}
\medskip
The following two lemmas are essential for the proof of Theorem 2.
\medskip
{\bf Lemma 3.} {\it Let $\,\kappa\geq\aleph_1\,$ be a 
regular cardinal. For $\,\theta\in\Omega\,$ consider the space 
$\;X=[0,\theta]\,$. If $\,\gamma\in C_\kappa(X)\,$
then either $\;\Omega_\kappa(\gamma)=[0,o(\kappa)[\;$
or $\;\Omega_\kappa(\gamma)=\{0\}\,$.}
\medskip
{\it Proof.} Since $\,\gamma\in C_\kappa(X)\,$, 
if $\,\alpha_1<\gamma\,$ and $\,\alpha_2\in\Omega\,$ 
and $\;|[\alpha_1,\alpha_2]|<\kappa\;$
then $\;[\alpha_1,\alpha_2]\subset [0,\gamma[\,$.
Clearly, 
if $\;\sup A\,\not=\,\gamma\;$ whenever $\;\emptyset\not=A\subset [0,\gamma[\;$
and $\,|A|<\kappa\,$ then $\;\Omega_\kappa(\gamma)=\{0\}\,$.
So assume that there is a nonempty set $\,A\subset [0,\gamma[\;$ such that
$\,|A|<\kappa\,$ and $\;\sup A\,=\,\gamma\,$.
For $\,\xi\in\Omega\,$
put $\;U_\xi\,:=\,\bigcup_{\alpha\in A}[\alpha,\alpha+\omega^\xi]\,$.
If $\,\xi<o(\kappa)\,$ then
$\;|[\alpha,\alpha+\omega^\xi]|=|[0,\omega^\xi]|<\kappa\;$
for every $\,\alpha\in A\,$ and hence
$\;U_\xi\subset [0,\gamma[\;$ 
and hence $\;\sup U_\xi\,=\,\gamma\,$.
Thus $\,|U_\xi|<\kappa\,$ and (by Lemma 2) $\;\gamma\in U_\xi^{(\xi)}\;$
for every $\,\xi<o(\kappa)\,$ and hence $\;\Omega_\kappa(\gamma)
=[0,o(\kappa)[\,$, {\it q.e.d.}
\medskip
{\it Remark.} In Lemma 3 the case 
$\;{\Psi}[[0,\theta],\kappa]=\{0,o(\kappa)\}\;$ may occur, for example if 
$\,\kappa=\aleph_1\,$ and $\,\theta=\omega_1\cdot\omega\,$.
\medskip\smallskip
{\bf Lemma 4.} {\it Let $\,(X,\prec)\,$ be a linearly ordered set
equipped with the order topology
and assume that the space $\,X\,$ is scattered. 
Let $\,0\not=\xi\in \Omega\,$ 
and let $\,\kappa\,$ 
be a regular cardinal number with $\;\kappa>|\omega^\xi|\,$.
Let $\;x,y,z\;$ be three points
in $\,X\,$ with $\;x\prec z\prec y\;$ so that 
$\;[x,z]_\prec\,$ 
is an order-isomorphic copy of $\,[0,\omega^\xi]\,$
and $\;[z,y]_\prec\,$ 
is an order-isomorphic copy of 
$\;[0,o(\kappa)]^*\,.$
Then $\;C_\kappa(X)\cap\,]x,y[_\prec\;=\,\{z\}\;$ 
and $\;\sup\,\Omega_\kappa(z)\,=\,\xi\,$.}
\medskip
{\it Proof.} Clearly, $\,z\,$ is
the only $\kappa$-condensation point of $\,X\,$ 
strictly between $\,x\,$ and $\,y\,$. 
Since $\,\kappa\,$ is regular, 
there is no set $\,A\subset [z,y]_\prec\,$
with $\,|A|<\kappa\,$ and $\,z\in A'\,$.
Therefore, if $\,0\not=\alpha\in\Omega\,$ and 
$\;z\in A^{(\alpha)}\;$
for a point set $\,A\,$ in the space $\,X\,$
with $\;|A|<\kappa\;$
then we already have $\;z\in (A\cap[x,z]_\prec)^{(\alpha)}\,$.
On the other hand, 
$\;([x,z]_\prec)^{(\xi)}=\{z\}\;$ by Lemma 2
and $\;|[x,z]|=|[0,\omega^\xi]|<\kappa\,$.
Consequently, $\;\Omega_\kappa(z)\,=\,[0,\xi]\;$
and hence  $\;\sup\,\Omega_\kappa(z)\,=\,\xi\,$, {\it q.e.d.}
\bigskip
Now we are ready to prove Theorem 2.
Let $\,\kappa\,$ be a cardinal with 
$\,\kappa>\aleph_0\,$
and put $\;K=[\omega,o(\kappa)]\,$.
Let $\,{\cal G}\,$ be the family 
of all nonempty sets $\,S\,$ of successor ordinals $\,\alpha+1\,$
with $\;\alpha\,\in\,K\cap{\cal L}\!\setminus\!\{o(\kappa)\}\,$.
So if $\,\xi\in S\in {\cal G}\,$ then $\;|\xi|=|[0,\omega^\xi]|<\kappa\,$.
Clearly, $\;|{\cal G}|=2^\kappa\,$.
For every set $\,S\in {\cal G}\,$ let
\medskip
\centerline{$H_S\;:=\;K\times\{0\}
\;\cup\;
\bigcup\limits_{\xi\in S}(\{\xi\}\!\times\![0,\omega^\xi]\;   
\cup\,\{\xi+1\}\!\times\! [0,o(\kappa)[^*)$}
\smallskip
and
\smallskip
\centerline{$G_S\;:=\;K\times\{0\}
\;\cup\;
\bigcup\limits_{\xi\in S}(\{\xi\}\!\times\! [0,\omega^\xi]\;\cup\,
\{\xi+1\}\!\times\! [0,o(|\xi|^+)[^*)$}
\medskip
be equipped with the lexicographic ordering.
Illustratively, the linearly ordered set $\,H_S\,$ resp.~$\,G_S\,$
is constructed from the well-ordered set 
$\,K\,$ by replacing $\,\xi\,$
with a copy of $\,[0,\omega^\xi]\;$ 
and $\,\xi+1\,$ with a copy
of $\,[0,o(\kappa)[^*\,$ resp.~$\,[0,o(|\xi|^+)[^*\,$
for each $\,\xi\in S\,$. 
\medskip
Then the corresponding linearly ordered spaces $\,H_S\,$  and $\,G_S\,$
are of size $\,\kappa\,$ and it is evident that 
all these spaces are scattered.
They are also compact since the ordering 
is complete with a maximum and a minimum (cf.~[10] 39.7).
We claim that the spaces $\;H_S\,(S\in{\cal G})\;$
are mutually non-homeomorphic if $\,\kappa\,$ is regular
and that the spaces $\;G_S\,(S\in{\cal G})\;$
are mutually non-homeomorphic if $\,\kappa\,$ is singular.
\medskip\smallskip
Assume firstly that $\,\kappa\,$ is regular
and let $\,S\in {\cal G}\,$ and consider the space $\,H_S\,$.
Clearly we have $\;(o(\kappa),0)\in C_\kappa(H_S)\;$ and
$\;\Omega_\kappa((o(\kappa),0))=\{0\}\,$. 
Obviously, $\;(\gamma,0)\in C_\kappa(H_S)\;$
if and only if $\,\gamma=o(\kappa)\,$
or $\;\gamma\,=\,\sup(S\cap[0,\gamma[)\;$
where $\;S\cap[0,\gamma[\,\not=\emptyset\,$.
If $\,(o(\kappa),0)\not=(\gamma,0)\in C_\kappa(H_S)\,$
then $\;\Omega_\kappa((\gamma,0))=[0,o(\kappa)[\;$ and hence
$\;\sup \Omega_\kappa((\gamma,0))=o(\kappa)\;$
because if $\;\xi\in S\cap[0,\gamma[\;$
then  $\;\{\xi+1\}\!\times\![0,\omega^\alpha]^*\,\subset\,
\{\xi+1\}\!\times\![0,o(\kappa)[^*\;$
and $\;|\bigcup\{\,\{\xi\!+\!1\}\!\times\![0,\omega^\alpha]^*\;|\;
\xi\in S\cap[0,\gamma[\,\}|<\kappa\;$
for arbitrarily large exponents $\,\alpha<o(\kappa)\,$.
In view of Lemma~4 we have
$\;C_\kappa(H_S)\setminus K\!\times\!\{0\}\,=\,
\{\,(\xi,\omega^\xi)\;|\;\xi\in S\,\}\;$
and $\;\sup\,\Omega_\kappa((\xi,\omega^\xi))\,=\,\xi\;$
for every $\,\xi\in S\,$. Therefore,
\medskip
\centerline{$\;S\,=\,{\Psi}[H_S,\kappa]\setminus\{0,o(\kappa)\}\;$}
\medskip
for every $\,S\in {\cal G}\,$
and this settles Theorem 2 for {\it regular} $\,\kappa>\aleph_0\,$.
\medskip\smallskip
Assume now that $\,\kappa\,$ is a singular cardinal number
and let $\,{\cal R}\,$ denote the set of all regular uncountable cardinals
smaller than $\,\kappa\,$. 
We claim that
for every $\,S\in{\cal G}\,$ we have 
\medskip\smallskip
\centerline{$S\,\;=\,\;\big(\bigcup\limits_{\lambda\in {\cal R}}
{\Psi}[G_S,\lambda]\,\big)
\setminus{\cal L}\;=\,\;\big(\bigcup\limits_{\lambda\in {\cal R}}
{\Psi}[G_S,\lambda]\,\big)
\setminus(\{0\}\cup\{\,o(\lambda)\;|\;\lambda\in{\cal R}\,\})\;.$} 
\medskip
On the one hand, if $\,\xi\in S\,$ 
then $\;\xi\not\in{\cal L}\;$ and $\;|\xi|^+\in {\cal R}\;$ and
$\,(\xi,\omega^{\xi})\,$ 
is a $\,|\xi|^+$-condensation point with
$\;\sup\,\Omega_{|\xi|^+}((\xi,\omega^\xi))\,=\,\xi\;$
in view of Lemma 4.
\medskip
On the other hand, let $\,y\,$
be a $\lambda$-condensation 
point in $\,G_S\,$ 
where $\,\lambda\in {\cal R}\,$
and assume firstly that $\;y\,\not\in\,K\times\{0\}\,$.
Then $\,y\,$ lies in 
$\;B_\xi\,:=\,\{\xi\}\times\,[0,\omega^\xi]\,\cup\,
\{\xi+1\}\times[0,|\xi|^+[^*\;$ for some $\,\xi\in S\,$.
Since the points $\;\min B_\xi\,=\,(\xi,0)\;$ and
$\;\max B_\xi\,=\,(\xi+1,0)\;$ 
are isolated in the space $\,G_S\,$,
the point $\,y\,$ must be a $\lambda$-condensation 
point in the space $\,B_\xi\,$, 
whence $\,\lambda\leq |B_\xi|=|\xi|^+\,$. 
\smallskip
In the case that $\,\lambda=|\xi|^+\,$
we must have $\;y=(\xi,\omega^\xi)\;$
and hence $\;\sup \Omega_\lambda(y)=\xi\in S\;$ by Lemma~4.
In the case that $\,\lambda<|\xi|^+\,$
the point $\,y\,$ must be the maximum resp.~minimum 
of an isomorphic copy of $\;[0,\gamma]\,$ resp.~$\;[0,\gamma]^*\,$
within the linearly ordered set $\,B_\xi\,$
where $\,\gamma\,$ is a $\lambda$-condensation point in the 
space $\,[0,\gamma]\,$, whence 
$\;\sup \Omega_\lambda(y)\in\{0,o(\lambda)\}\;$ by Lemma 3.
\smallskip
Assume secondly that $\,y=(x,0)\,$ for $\,x\in K\,$.
If $\,x\,$ is a $\lambda$-condensation point in the basic
space $\,K\,$ then it is clear that 
in the space $\,G_S\,$ we also have 
$\;\sup \Omega_\lambda(y)\in\{0,o(\lambda)\}\,$.
If $\;x\not\in C_\lambda(K)\;$ 
then $\,y\in C_\lambda(G_S)\,$ 
forces $\,x\,$ to be the supremum 
of a set 
\smallskip
\centerline{$\;\tilde S\,\subset\,
\{\,\xi\in S\;|\;\xi<x\;\land\;|\xi|^+=\lambda\,\}\;$}
\smallskip
with $\,|\tilde S|<\lambda\,$ and therefore 
(by the same argument as for the space $\,H_S\,$) 
we must have $\,\Omega_\lambda(y)=[0,o(\lambda)[\;$
and hence $\,\sup \Omega_\lambda(y)=o(\lambda)\,$.
\smallskip
So in any case the ordinal number $\,\sup \Omega_\lambda(y)\,$ lies 
in $\;S\cup\{0,o(\lambda)\}\;$ if 
$\;y\in C_\lambda(G_S)\;$ for $\,\lambda\in{\cal R}\,$.
This concludes the proof of  Theorem 2.
\bp\mp
{\bff 9. Completions and compactifications}
\mp
In this short, final section we present two nice applications
of Theorem 1 and Theorem 2.
\mp
If $\,X\,$ is a scattered Hausdorff space then it is plain 
that the set $\,X\setminus X'\,$ of all isolated points is 
dense. Moreover, $\,X\setminus X'\,$ is the intersection of all
dense subsets of $\,X\,$. In particular, $\,|X\setminus X'|\,$ 
is the {\it density} of the scattered space $\,X\,$.
Consequently, if the scattered space $\,X\,$ is metrizable
then $\,|X\setminus X'|=|X|\,$ because {\it weight} and {\it density}
of a metric space are always identical.
Therefore, from Theorem 1 we derive the following enumeration result
about completions of discrete metric spaces. 
\mp
{\bf Corollary 3.} {\it The topology of an infinite discrete 
space $\,S\,$ can be generated by 
$\,2^{|S|}\,$ ultra\-metrics $\,d\,$ such that 
the completions of the metric spaces $\,(S,d)\,$ are 
mutually non-homeomorphic scattered, ultrametric (linearly ordered)
spaces (of weight and size $\,|S|\,$).}
\mp
{\it Remark.} Size and weight of any scattered completion 
of a discrete metric space $\,S\,$ 
must coincide with $\,|S|\,$ since size and density of
a scattered metric space are always identical.        
\mp\sp
Similarly, from the proof of Theorem 2 we can derive 
an enumeration theorem about compactifications of discrete spaces.
While $\;|X\setminus X'|<|X|\;$ is possible for compact and
scattered Hausdorff spaces $\,X\,$ (consider for example 
the one point compactification of space {\bf 65} in [10]),
in the proof of Theorem 2 it is evident that 
$\;|X\setminus X'|=|X|=\kappa\;$
whenever $\;X\in\{H_S, G_S\}\;$ for $\,S\in{\cal G}\,$.
This is clearly enough to settle the following enumeration result.
\mp
{\bf Corollary 4.} {\it Every uncountable discrete 
space $\,S\,$ has precisely $\,2^{|S|}\,$ scattered and linearly 
ordered compactifications of size (and weight) $\,|S|\,$ up to homeomorphism.}
\mp
{\it Remark.} As already pointed out, the statement in Corollary 4 would be
unprovable for a countably infinite discrete space $\,S\,$.
However, in view of [9] it is clear that any 
countably infinite discrete space $\,S\,$ has 
precisely $\,\aleph_1\,$ countable compactifications up to homeo\-morphism,
and they all are scattered and linearly 
ordered spaces. Furthermore, it is worth mentioning 
that any countably infinite discrete space $\,S\,$ has 
$\,2^{\aleph_0}\,$ mutually non-homeomorphic 
{\it uncountable} compactifications which all are also
{\it linearly ordered} and {\it metrizable} spaces (see [8] Theorem 7).  
\bp\bp
{\bff References}
\medskip
[1] Engelking, R.,  {\it General Topology}, revised and completed edition.
Heldermann 1989. 
\smallskip
[2] Hodel, R.E., {\it The number of metrizable spaces}, 
Fund.~Math.~{\bf 115} (1983),
127-141.
\smallskip
[3] Kechris, A., {\it Classical Descriptive Set Theory}, Springer 1995.
\sp
[4] Kuba, G.: {\it Counting topologies.}
Elemente d.~Math. {\bf 66} (2011), 56-62.
\sp
[5] Kuba, G., {\it Counting metric spaces},
Archiv d.~Math.~{\bf 97} (2011), 569-578.
\smallskip
[6] Kuba, G., {\it Counting linearly ordered spaces},
Coll.~Math.~{\bf 135} (2014), 1-14.
\sp
[7] Kuba, G., {\it Counting ultrametric spaces},
Coll.~Math.~{\bf 152} (2018), 217-234.
\sp
[8] Kuba, G., {\it On the variety of Euclidean point sets},
arXiv:2004.11101v1 [math.GN] (2020).
\sp
[9] Mazurkiewicz, S., and  Sierpi\'nski, W., {\it Contribution … la topologie
des ensembles}

\rl{{\it  d'nom\-brables}, Fund.~Math.~{\bf 1} (1920), 17-27.}

[10] Steen, L.A., and Seebach Jr., J.A.,
{\it Counterexamples in Topology,} 
Dover 1995.
\bp\mp
{\sl Author's address:} Institute of Mathematics. 

University of Natural Resources and Life Sciences, Vienna, Austria. 

{\sl E-mail:} {\tt gerald.kuba@boku.ac.at}

\bigskip\medskip
\hrule
\bigskip\bigskip
{\bf Theorem 1 and its proof is contained in the author's paper [7].}
\medskip
{\bf Theorem 2 and its proof is contained in the author's paper [6].}
\bye
\end